\documentclass[11pt,
]{article}

\usepackage{amssymb,amstext,amsthm,latexsym}
\usepackage{amsmath}
\usepackage{amsfonts}
\usepackage{mathrsfs}
\usepackage{mparhack}
\usepackage[bbgreekl]{mathbbol}

\newcommand{\nek}{\newcommand}
\nek{\renek}{\renewcommand}

\makeatletter
\if@twoside%
\else\fi%
\makeatother

\nek{\punk}{\subsection}
\nek{\parf}{\subsection}

\nek{\skr}{\mathscr}
\nek{\cP}{\skr P}

\nek{\vyk} [1] {}

\nek{\imar}[1]{} 
\vyk{
{\marginpar[
\flushright\footnotesize%
$\mtho\longrightarrow$\\%
\vspace{-1ex}{#1}\vspace*{1ex}]%
{
\flushleft\footnotesize%
$\mtho\longleftarrow$\\%
\vspace{-1ex}{#1}\vspace*{1ex}}%
}%
}

\nek{\imae} {\imar}

\nek{\itsep}{\itemsep=0.3ex plus 0.1ex minus 0.1ex}
\nek{\tenu}[1]{
\def\theenumi{#1}
\def\labelenumi{\theenumi}\itsep
}

\theoremstyle{plain}

\newtheorem{theore}             {Theorem} [subsection] 
\newtheorem{corollar}  [theore]{Corollary}
\newtheorem{propo}     [theore]{Proposition}
\newtheorem{lemm}      [theore]{Lemma}
\newtheorem{lemt}      {Lemma}[theore]
\newtheorem{cla}       [theore]{Claim} 
\newtheorem{clt}       {Claim} [theore]

\theoremstyle{definition}
\newtheorem{defn}      [theore]{Definition}
\newtheorem{rem} [theore]   {Remark} 
\newtheorem*{solu}{Solution}      
\newtheorem*{probl}{Problem}      
\newtheorem{exz} {Example} 
\newtheorem*{prF}{Proof}               

\nek{\thsp}{\hspace{0.1ex plus \mathsurround}}

\nek{\bex}{\begin{exz}}
\nek{\eex}{\end{exz}}
\nek{\bprl}{\begin{probl}}
\nek{\eprl}{\end{probl}}
\nek{\bsol}{\begin{solu}}
\nek{\esol}{\end{solu}}
\nek{\bpro}{\begin{propo}}
\nek{\epro}{\end{propo}}
\nek{\bcor}{\begin{corollar}}
\nek{\ecor}{\end{corollar}}
\nek{\bdf} {\begin{defn}} 
\nek{\eDf} {\end{defn}}
\nek{\edf} {\end{defn}}
\nek{\edF} {\end{defn}}
\nek{\ble} {\begin{lemm}}
\nek{\ele} {\end{lemm}}
\nek{\blt} {\begin{lemt}}
\nek{\elt} {\end{lemt}}
\nek{\bre}{\begin{rem}}
\nek{\ere}{\end{rem}}
\nek{\bte} {\begin{theore}}
\nek{\ete} {\end{theore}}
\nek{\bpf} {\begin{prF}} 
\nek{\epf} {\qed\end{prF}} 
\nek{\ePf} {\end{prF}} 
\nek{\qeD} {\qed}
\nek{\qeDD} [1] 
{\hfill\hbox{\qed~({\small #1\/}\hspace{0.1ex})}}
\nek{\epF} [1] {\qeDD{#1}\end{prF}} 

\nek{\ben}{\begin{enumerate}\itsep}
\nek{\een}{\end{enumerate}}
\nek{\bit}{\begin{itemize}\itsep}
\nek{\eit}{\end{itemize}}
\nek{\bde}{\begin{description}\itsep}
\nek{\ede}{\end{description}}
\nek{\bay}{\begin{array}}
\nek{\eay}{\end{array}}
\nek{\bce}{\begin{center}}
\nek{\ece}{\end{center}}
\nek{\bqo}{\begin{quotation}\noi}
\nek{\eqo}{\end{quotation}}
\nek{\ZFC}{{\bf ZFC}}
\nek{\ZF}{{\bf ZF}}

\nek{\iesp}{\hspace{0.3ex}}
\nek{\resp}{\hspace{0.25ex}}
\nek{\ie} {\hbox{\sl i.\iesp e.}}
\nek{\eg} {\hbox{\sl e.\iesp g.}}
\nek{\ea} {\hbox{\sl e.\iesp a.}}
\nek{\noo}
{\hbox{\sl w.\iesp l.\iesp o.\iesp g.}}
\nek{\vrt} {\hbox{w.\iesp r.\iesp t.}}
\nek{\ddd}[1]{$\mtho\hspace{0.2ex}{#1}$-\hspace{0.0ex}}
\nek{\dd}{\ddd}
\nek{\ran}  {\mathop{\tt ran}}
\nek{\otp} [1] {\mathop{\rm otp}(#1)}
\nek{\dom}  {\mathop{\tt dom}}
\nek{\tsup} {\mathop{\tt sup}}
\nek{\tmax} {\mathop{\tt max}}
\nek{\tlim} {\mathop{\tt lim\hspace{0.3ex}}}
\nek{\tmin} {\mathop{\tt min}}
\nek{\al} {\alpha}
\nek{\La}{\Lambda}
\nek{\la}{\lambda}
\nek{\sg} {\sigma}
\nek{\da} {\delta}
\nek{\Sg} {\Sigma}
\nek{\ve}{\varepsilon}
\nek{\vpi}{\varphi}
\nek{\vt} {\vartheta}
\nek{\om} {\omega}
\nek{\Om} {\Omega}
\nek{\lom}{^{<\om}}
\nek{\lomi}{^{<\omi}}
\nek{\za} {\zeta}
\nek{\omi} {\om_1}
\nek{\ali} {\aleph_1}
\nek{\alo} {{\aleph_0}}

\nek{\fs}[2]{{\hspace*{0.3ex}\boldsymbol\Sigma}^{#1}_{#2}}
\nek{\fp}[2]{{\boldsymbol\Pi}^{#1}_{#2}}
\nek{\fd}[2]{{\boldsymbol\Delta}^{#1}_{#2}}

\nek{\id}[2]{{\varDelta}^{#1}_{#2}}
\nek{\ip}[2]{{\varPi}^{#1}_{#2}}
\nek{\is}[2]{{\varSigma}^{#1}_{#2}}
\nek{\iSg}{\varSigma}
\nek{\iPi}{\varPi}

\nek{\BBB}{\hspace{0.01ex}}
\nek{\dR}{{\BBB{\mathbb R}\BBB}}
\nek{\dP}{{\BBB{\mathbb P}\BBB}}

\nek{\bn}{\om^\om}
\nek{\sus} {\mathopen{\exists\hspace{0.35ex}}}
\nek{\kaz} {\mathopen{\forall\hspace{0.35ex}}}
\nek{\imp} {\Longrightarrow} 
\nek{\mpi} {\Longleftarrow} 
\nek{\eqv} {\Longleftrightarrow} 
\nek{\leqv} {\;\eqv\;} 
\nek{\limp} {\;\imp\;} 
\nek{\ti}  {\times} 
\nek{\sq}  {\subseteq}
\nek{\su}  {\subset}
\nek{\sneq}{\subsetneqq}
\nek{\we}  {{\mathbin{\hspace{0.15ex}^\wedge}}}
\nek{\obr} {^{-1}}
\nek{\dif} {\smallsetminus}
\nek{\res} {\mathbin{\restriction}}
\nek{\lef} {\preccurlyeq}
\nek{\gef} {\succcurlyeq}
\nek{\pu}  {\varnothing}
\nek{\iy}  {\infty}
\nek{\piy} {+\iy}
\nek{\nin} {\not\in}
\nek{\onto}{\stackrel{\text{\rm onto}}{\longrightarrow}}
\nek{\ang} [1] {\langle #1\rangle}
\nek{\stk} [2] {\ang{#1\hspace{0.3ex};\hspace{0.1ex}#2}}
\nek{\sis} [2] {\ans{#1}_{#2}}
\nek{\ans} [1] {\{\hspace{0.01ex}#1\hspace{0.01ex}\}}

\nek{\zz} {\linebreak[0]} 
\nek{\ens} [2] {\ans{{#1\hspace{0.5ex}{:}}\zz\hspace{0.5ex}#2}}
\nek{\itla} {\item\label}

\nek{\ubf}{\fontseries{b}\selectfont}
\nek{\TS}{\textstyle}
\nek{\yo} {,\linebreak[0]}
\nek{\yi} {\hspace{0.2ex},\linebreak[0]\hspace{0.2ex}}
\nek{\yd} {\hspace{0.2ex},\linebreak[0]\:}
\nek{\yt} {\hspace*{0.2ex},\linebreak[0]\;} 
\nek{\lex} {<_{\text{\tt lex}}}
\nek{\lexe} {\leqslant_{\text{\tt lex}}}

\nek{\snos} [1] {\,\footnote{\hspace*{2pt}#1}}
\nek{\renu}{\tenu{{\rm(\roman{enumi})}}}
\nek{\fenu}{\tenu{{\rm(\fnsymbol{enumi})}}}
\nek{\Renu}{\tenu{{\rm(\Roman{enumi})}}}
\nek{\rit} [1] {{\it#1\/}}
\nek{\lam} [1] {\label{#1}\hspace*{-3pt}\imar{#1}}%
\nek{\las} [1] {\label{#1}\imar{#1}}%

\nek{\atc}{\addtocounter{enumi}{1}}
\nek{\wo} {\text{\ubf WO}}
\nek{\bez} {\dif}
\nek{\AC}   {{\text{\bf AC}}}
\nek{\DC}   {{\text{\bf DC}}}
\nek{\ROD}  {{\text{\bf{ROD}}}}
\nek{\hrod}{{\text{\bf HROD}}}
\nek{\OD} {\mathbf P}

\nek{\od} {{\text{\bf OD}}}

\nek{\rL} {\text{\ubf L}}
\nek{\rV} {\text{\ubf V}}
\nek{\rU} {\text{\ubf U}}

\nek{\led} {\leq^\ast}
\nek{\ld} {<^\ast}

\nek{\pqo} {PQO}

\nek{\Ord}  {\mathop{\tt Ord}}
\nek{\supp}  {\mathop{\tt supp}}
\nek{\card}  {\mathop{\tt card}}

\nek{\ba} {\beta}
\nek{\ga} {\gamma}
\nek{\ka} {\kappa}

\nek{\np}{\newpage}

\nek{\mtho}{\mathsurround=0mm}
\nek{\msur}{\hspace{-1\mathsurround}}
\nek{\noi}{\noindent}
\nek{\vom}{\vspace{1mm}}
\nek{\vim}{\vspace{-1mm}}
\nek{\vtm}{\vspace{2mm}}

\nek{\eqr} {equivalence relation}

\nek{\rE} {\mathrel{\mathsf E}}

\nek{\ek}[2] {[#1]_{{#2}}}

\nek{\eke}[1] {\ek{#1}{\rE}}

\nek{\qand}{\quad\text{and}\quad}

\nek{\dX}{{\BBB{\mathbb X}\BBB}}
\nek{\dN}{{\BBB{\mathbb N}\BBB}}
\nek{\dB}{{\BBB{\mathbb B}\BBB}}
\nek{\dF}{{\BBB{\mathbb F}\BBB}}
\nek{\fun}{{\BBB{\mathbb {Fun}}\BBB}}

\nek{\doP}  [1] {{#1}^\complement}

\nek{\curle}{\preccurlyeq}
\nek{\cle}{\curle}
\nek{\cl} {\prec}
\nek{\ncle}{\not\cle}
\nek{\ncl}{\not\cl}

\nek{\gh}{\mathbb P}
\nek{\ghd}{\gh^2}

\nek{\fdt}{\hbox{\raisebox{-0.25ex}{\LARGE\bf.}}}
\nek{\bdot}[1] 
{\raisebox{-0.07ex}{\mtho$\stackrel{\fdt}{#1}$}}

\nek{\dox}{\bdot{\text{\bf x}}}
\renek{\dox}{\bdot{\boldsymbol x}}
\nek{\doxl}{\dox_{\tt le}}
\nek{\doxr}{\dox_{\tt ri}}

\nek{\dn}{2^\om}
\nek{\dpe} {\mathord{\drof\gh\rE}}
\nek{\dpew}{\mathord{\drow{\gh}\rE}}

\nek{\trof} [3] {{#1}\ti_{#2}{#3}}

\nek{\drof} [2] {{#1}\ti_{#2}{#1}}
\nek{\drow} [2] {{#1}\ti_{#2}^{\text{weak}}{#1}}

\nek{\ck} {\om_1^{\text{\sc ck}}}

\nek{\Eo}  {\rE_{\text{\sf0}}}
\nek{\Fo}  {\rF_{\text{\sf0}}}
\nek{\nE}  {\mathbin{{\not\hspace{-0.35ex}\sf E}}}

\nek{\bcl} {\begin{cla}}
\nek{\ecl} {\end{cla}}
\nek{\bct} {\begin{clt}}
\nek{\ect} {\end{clt}}

\nek{\bV}{{\mathbf V}}
\nek{\bL}{{\bf L}}

\nek{\gM} {\mathfrak M}
\nek{\gP} {\mathfrak P}
\nek{\cM} {\skr M}

\nek{\ccs} {}
\nek{\cF}{{\ccs{\skr F}\ccs}}
\nek{\cS}{{\ccs{\skr S}\ccs}}
\nek{\cD}{{\ccs{\skr D}\ccs}}
\nek{\cN}{{\ccs{\skr N}\ccs}}
\nek{\cO}{{\ccs{\skr O}\ccs}}

\nek{\cf} [1] {\cF_{#1}}
\nek{\cfx} [1] {\cF_#1}
\nek{\cfd} [2] {\cF_{#1}(#2)}

\nek{\etc} {{\sl etc}}

\nek{\bus}{\begin{equation}}   
\nek{\eus}{\end{equation}}

\nek{\pp} [2] {\gh^{#1}_{#2}}

\nek{\dpd} [1] {\gh\ti_{#1}\gh}
\nek{\dpx} [1] {\gh\ti_{\rE_{#1}}\gh}

\nek{\dpw} {\gh(W)}
\nek{\dpwe} {\gh(W)\ti_{\rE}\gh(W)}

\nek{\ups} {\varUpsilon}

\nek{\apr} {\approx}
\nek{\napr}{\not\apr}

\nek{\gp}{\mathfrak p}
\renek{\gp}{\mathbb p}

\nek{\aenu}{\tenu{{\rm(\arabic{enumi})}}}
\nek{\Aenu}{\tenu{{\rm(\Alph{enumi})}}}
 
\nek{\doxlp}{\dox{}'_{\tt le}}
\nek{\doxrp}{\dox{}'_{\tt ri}}

\nek{\bfit}{\bfseries\itshape}

\nek{\esn}{\rS_{\ans{1/n}}}
\nek{\rS}  {\mathbin{\sf S}}

\nek{\wh}{\widehat}
\nek{\wY}{\widehat Y}
\renek{\wY}{C}

\nek{\lr}{LR}
\nek{\rl}{RL}

\nek{\vx} {\vec x}

\nek{\Xa}{X^*}
\nek{\Ua}{U^*}
\nek{\Uo}{U_0}
\nek{\Va}{V^*}
\nek{\Do}{D_0}

\nek{\rEF} {\mathrel{\mathsf E_\cF}}
\nek{\refx} [1] {\mathrel{\mathsf E_{\cfx#1}}}

\nek{\osm} {\dd\Omega SM}

\nek{\col} [1] {\text{\ubf Coll}(\om,#1)}
\nek{\coll} [1] {\col{{<\hspace*{0.1ex}}#1}}

\nek{\dodg} {\od-generic}
\nek{\pge} {\dd\OD generic}
\nek{\ege} {\dd{(\spe)}generic}
\nek{\ega} [1] {\dd{(\spa{#1})}generic}

\nek{\spe} {\mathord{\drof\OD\rE}}
\nek{\spa} [1] {\mathord{\drof\OD{\rE_{#1}}}}

\nek{\spei} 
{\mathord{\drof{\odi\hspace*{-0.3ex}}\rE\hspace*{0.3ex}}}

\nek{\spai} [1] 
{\mathord{\drof{\odi\hspace*{-0.3ex}}{\rE_{#1}}\hspace*{0.3ex}}}

\nek{\xle} {x\ile}
\nek{\xri} {x\iri}

\nek{\odf} {\od-force}
\nek{\pfo} {\dd\OD force}
\nek{\efo} {\dd{(\spe)}force}
\nek{\efa} [1] {\dd{(\spa{#1})}force}

\nek{\odsw} [1] {\OD_{\sq#1}}

\nek{\ile} {_{\text{\tt le}}}
\nek{\iri} {_{\text{\tt ri}}}

\nek{\odk} {\dd\od1st-countable\/}

\nek{\pwod}  [1] {\cP_{\text{\tt OD}}(#1)}
\nek{\pws}  [1] {\cP(#1)}
\nek{\odi} {\OD^*}
\nek{\koh} {\text{\sc Coh}}

\nek{\dotx}{\check x}
\nek{\dotY}{\check Y}

\nek{\ibn} [1] {\skr N_{#1}}

\nek{\cenu}{\tenud{{%
$\mtho\arabic{enumi}^\circ$%
}}}

\nek{\odw} {\odsw{\ppw}}
\nek{\odwp} {\odsw{\ppl'}}
\nek{\odwx} {\odsw{B}}
\nek{\odwb} {\odsw{\ppl}}
\nek{\odwe} {\odw\ti_{\rE}\odw}
\nek{\odwep} {\odwp\ti_{\rE}\odwp}
\nek{\odwex} {\odwx\ti_{\rE}\odwx}

\nek{\odwbi} 
{\mathord{\odi_{\sq\ppl}}}

\nek{\tenud}[1]{
\def\theenumi{#1}
\def\labelenumi{\theenumi.}\itsep
}

\nek{\bse} {2\lom}
\nek{\ilom}{\omi\lom}

\nek{\vT} {\Theta}

\nek{\Lom} {^{<\Om}} 
\nek{\msl} {\mathbin{<_{\text{lex}}}} 

\nek{\cof}  {\mathop{\tt cof}}

\renek{\od} {\text{OD}}
\nek{\rod} {\text{ROD}}

\nek{\wX}{C}
\nek{\mek}{\cle}
\nek{\eee}{\apr}
\nek{\omck}{\omi^{CK}}
\nek{\tfu}{\cF}
\nek{\hop} {HOP}
\nek{\strk}{\stk}

\nek{\lra}{\to}

\nek{\mel}{\lexe}

\nek{\eqf} {\equiv}
\nek{\cj}{\land}
\nek{\dm}{$$}

\nek{\prx}{{\text{\tt pr}}_x\hspace{1pt}}
\nek{\pry}{{\text{\tt pr}}_y\hspace{1pt}}
\nek{\pri} {\dom} 
\nek{\prt} {\ran} 

\nek{\emps}{\pu}
\nek{\Xin}  {X_\infty}
\nek{\kV}  {\rV}
\nek{\kvp}{\kV^+}

\nek{\dpp}{\dP^+_2}
\nek{\dpm}{\dP^-_2}
\nek{\dpt}{\dP^2_{\eqf}}

\nek{\Ups} {\ups}

\nek{\nEo} {\mathrel{\not{{\hspace{-0.4ex}\rE}}_0}}
\nek{\meo} {\mathrel{\leq_0}}

\nek{\pone}{\hspace{-0.4ex}+\hspace{-0.4ex}1}

\nek{\ppl} {\boldsymbol B}
\nek{\ppw} {W}
\nek{\pmi} {\boldsymbol C}

\begin{document}

\title
{Linearization of partial quasi-orderings in the Solovay model 
revisited.
\thanks{Partial support of 
RFFI grant 13-01-00006 acknowledged.}}

\author 
{Vladimir~Kanovei\thanks{IITP RAS and MIIT,
  Moscow, Russia, \ {\tt kanovei@googlemail.com} --- contact author}  
\and
Vassily~Lyubetsky\thanks{IITP RAS,
  Moscow, Russia, \ {\tt lyubetsk@iitp.ru} 
}}

\date 
{\today}

\maketitle

\begin{abstract}
We modify arguments in \cite{ya:alin} to reprove a linearization 
theorem on real-ordinal definable 
partial quasi-orderings in the Solovay model.
\end{abstract}

\subsection{Introduction}
 
The following theorem is the main content of this note.

\bte
[in the Solovay model]
\lam{mt}
Let\/ $\cle$ be a\/ \ROD\ 
(real-ordinal definable) 
partial quasi-ordering on\/ $\bn$ and\/
$\apr$ be the associated \eqr.
Then exactly one of the 
following two conditions is satisfied$:$
\ben
\Renu
\itla{mt1} 
there is an antichain\/ $A\sq 2\lomi$ and a\/ $\ROD$ map\/ 
$F:\bn\to A$ such that \

$1)$ 
if\/ $a,b\in\bn$ then$:$ $x\cle y\imp F(x)\lexe F(y)$, \ and \ 

$2)$ 
if\/ $a,b\in\bn$ then$:$ ${x\napr y}\limp{F(x)\ne F(y)}\;;$ 

\itla{mt2} 
there exists a continuous $1-1$ map\/ 
$F:\dn\to\bn$ such that\/ 

$3)$ 
if\/ $a,b\in\dn$ then$:$ $a\meo b\imp F(a)\cle F(b)$, \ and \ 

$4)$ 
if\/ $a,b\in\dn$ then$:$ ${a\not\Eo b}\limp{F(a)\ncle F(b)}\;.$ 
\een
\ete 

Here $\lexe$ is the lexicographical order on sets of the 
form $2^\al\yi\al\in\Ord$ --- it linearly orders any antichain 
$A\sq 2\lomi,$ 
while $\meo$ is the partial quasi-ordering on $\dn$ 
defined so that $x\meo y$ 
iff $x\Eo y$ and either $x=y$ or $x(k)<y(k)$, where $k$ is 
the largest number with $x(k)\ne y(k)$.\snos
{Clearly
$\meo$ orders each \dd\Eo class similarly to 
the (positive and negative) integers, except for the class 
$\ek{\om\ti\ans0}{\Eo}$ ordered as $\om$ and
the class $\ek{\om\ti\ans1}{\Eo}$ ordered the inverse of $\om$.}

The proof of this theorem (Theorem 6) in \cite[Section 6]{ya:alin}) 
contains a reference to Theorem 5 on page 91 (top), 
which is in fact not immediately applicable in the Solovay model. 
The goal of this note is to present a direct and self-contained 
proof of Theorem~\ref{mt}. 

The combinatorial side of the proof follows the proof of a 
theorem on Borel linearization in \cite{k:blin}, in turn based 
on earlier results in \cite{hms,hkl}. 
This will lead us to \ref{mt1} in a weaker form, with a function 
$F$ mapping $\bn$ into $2^{\omega_2}.$  
To reduce this to an antichain in $2\lomi$, a compression lemma 
(Lemma~\ref{apal31} below)
is applied, which has no counterpart in the Borel case. 

Our general notation follows \cite{kanB,ksz}, but 
for the convenience of the reader, we add a review of notation.

%
\bde
\item[\rm PQO, \it partial quasi-order\/$:$] \ 
reflexive ($x\le x$) and transitive in the domain;


\item[\rm LQO,  
\it linear quasi-order\,$:$] \  
PQO and $x\le y\lor y\le x$  
in the domain; 

\item[\rm LO, \it linear order\,$:$]  \ 
LQO and $x\le y\land y\le x\imp x=y$; 

\item[\it associated equivalence relation\,$:$]  \ 
$x\apr y$ iff $x\le y\land y\le x$.

\item[\it associated strict ordering\,$:$]  \ 
$x<y$ iff $x\le y\land y\not\le x$;

\item[\it \lr\ (left--right) order preserving map$:$]
any map $f:\stk{X}{\le}\to\stk{X'}{\le'}$ such that 
we have $x\le y\imp f(x)\le' f(y)$ for all $x,y\in\dom f$;

\item[\it $\lex\yt\lexe\;:$]
the lexicographical LOs on sets of the form $2^\al\yd\al\in\Ord$, 
resp.\ strict and non-strict;

\item[\rm $\eke x=\ens{y\in\dom\rE}{x\rE y}$
(the \dd\rE\rit{class} of $x$) 
\ and \ 
$\eke X=\textstyle\bigcup_{x\in X}\eke x$]
--- \\ whenever  $\rE$ is an \eqr\ and $x\in\dom\rE$, $X\sq\dom\rE$.
\ede

\bre
\lam{ODonly}
We shall consider only the case of a parameterfree \od\ 
ordering $\cle$ in Theorem~\ref{mt}; the case of $\od(p)$ with 
a fixed real parameter $p$ does not differ much.
\ere

\subsection{The Solovay model and $\od$ forcing}
\las{smod}
                                                    
We start with a brief review of the Solovay model.
Let $\Om$ be an ordinal. 
Let \osm\ be the following hypothesis:
\bde
\item[\rm\osm:]
$\Om=\omi$, 
$\Om$ is strongly inaccessible in $\rL$, the 
constructible universe,   and
the whole universe $\rV$ is 
a generic extension of $\rL$ via the Levy collapse forcing 
$\coll\Om$, as in \cite{solMST}.
\ede
Assuming \osm, let $\OD$ be the set of all {\ubf non-empty} 
$\od$ sets $Y\sq \bn$.
We consider $\OD$ as a forcing notion 
(smaller sets are stronger).
A set $D\sq \OD$ is:
\bit
\item[$-$]\rit{dense}, \ 
iff for every $Y\in\OD$ there exists $Z\in D$, $Z\sq Y$;

\item[$-$]\rit{open dense}, \ 
iff in addition we have $Y\in D\imp X\in D$ whenever 
sets $Y\sq X$ belong to $\OD$;
\eit   
A set $G\sq \OD$ is {\pge}, \ iff \ 
1) 
if $X,Y\in G$ then there is a set $Z\in G$, $Z\sq X\cap Y$,   
\hspace*{0.3ex}and \ 
2) 
if $D\sq\OD$ is $\od$ and dense then $G\cap D\ne\pu$.

Given an $\od$ \eqr\ $\rE$ on $\bn,$ a 
\rit{reduced product} forcing 
notion $\spe$ consists of all sets of the form 
$X\ti Y,$ where $X\yi Y\in\OD$ and $\eke X\cap\eke Y\ne\pu$.
For instance $X\ti X$ belongs to $\spe$ whenever $X\in\OD$. 
The notions of sets dense and open dense in $\spe$, and 
\ege\ sets are similar to the case of $\OD$

A condition $X\ti Y$ in $\spe$ is \rit{saturated} iff 
$\eke X=\eke Y$.

\ble
\lam{sm6}
If\/ $X\ti Y$ is a condition in $\spe$ then there is a 
stronger saturated subcondition\/ $X'\ti Y'$ in $\spe$. 
\ele
\bpf
Let $X'=X\cap\eke Y$ and $Y'=Y\cap\eke X$. 
\epf

\bpro
[lemmas 14, 16 in \cite{ksol}] 
\lam{genx}
Assume \osm. 

If a set\/ $G\sq\OD$ is\/ \pge\ then the 
intersection\/ $\bigcap G=\ans{x[G]}$ consists of a 
single real\/ $x[G]$, called\/ \pge\ ---  
its name will be\/ $\dox$.

Given an $\od$ \eqr\ $\rE$ on $\bn,$ 
if a set\/ $G\sq\spe$ is\/ \ege\ then the 
intersection\/ $\bigcap G=\ans{\ang{\xle[G],\xri[G]}}$ 
consists of a single pair of reals\/ $\xle[G]\yi\xri[G]$, 
called an\/ \ege\ pair --- 
their names will be\/ $\doxl\yi\doxr$; 
either of\/ $\xle[G]\yi\xri[G]$ is separately\/ \pge.
\qed
\epro

As the set $\OD$ is definitely uncountable, 
the existence of \pge\ sets does not 
immediately follow from \osm\ by a cardinality argument. 
Yet fortunately $\OD$ is \rit{locally countable}, in a sense. 

\bdf
[assuming \osm]
\lam{odik}
A set $X\in\od$ is \rit{\odk} if the set 
$\pwod X=\pws X\cap\od$ of all $\od$ subsets of $X$ 
is at most countable. 
\edf

For instance, assuming \osm, the set $X=\bn\cap\od=\bn\cap\rL$ 
of all $\od$ reals is \odk. 
Indeed  
$\pwod X = \pws X\cap \rL$, and hence 
$\pwod X$ admits an $\od$  
bijection onto the ordinal $\om_2^{\rL}<\omi=\Om$.

\ble
[assuming \osm]
\lam{2co}
If a set $X\in\od$ is \rit{\odk} then the set\/ $\pwod X$ 
is\/ \odk either.
\ele
\bpf
There is an ordinal $\la<\omi=\Om$ and an $\od$ bijection 
$b:\la\onto\pwod X$. 
Any $\od$ set $Y\sq\la$ belongs to $\bL$, hence, the $\od$ 
power set $\pwod\la=\pws\la\cap\bL$ belongs 
to $\bL$ and $\card(\pwod\la)\le\la^+<\Om$ in $\bL$. 
We conclude that $\pwod\la$ is countable. 
It follows that $\pwod{\pwod X}$ is countable, as required.
\epf


\ble
[assuming \osm]
\lam{L*}
If\/ $\la<\Om$ then the set\/ $\koh_\la$ of all 
elements\/ $f\in\la^\om$, \dd{\col\la}generic over\/ $\rL$, 
is\/ \odk. 
\ele
\bpf
If $Y\sq\koh_\la$ is \od\ and $x\in Y$ then 
``$\dotx\in \dotY$'' is \dd{\col\la}forced over $\rL$. 
It follows that there is a set 
$S\sq\la\lom=\col\la\yt S\in\rL$, such 
that $Y=\koh_\la\cap\bigcup_{t\in S}\ibn t$, where 
$\ibn t=\ens{x\in\la\lom}{t\su x}$, 
a Baire interval in $\la\lom.$ 
But the collection of all such sets $S$ belongs to $\rL$ 
and has cardinality 
$\la^+$ in $\rL$, hence, is countable under \osm.
\epf

Let $\odi$ be the set of all \odk\ sets $X\in\OD$.
We also define 
$$
\spei=\ens{X\ti Y\in\spe}{X,Y\in \odi}.
$$ 

\ble
[assuming \osm]
\lam{den}
The set\/ $\odi$ is dense in\/ $\OD$, that is, 
if\/ $X\in\OD$ then there is a condition\/ $Y\in\odi$ 
such that\/ $Y\sq X$. 

If\/ $\rE$ is an\/ $\od$ \eqr\ on $\bn$ 
then the set\/ $\spei$ is dense in\/ $\spe$ 
and any\/ $X\ti Y$ in\/ $\spei$ is\/ \odk.
\ele
\bpf
Let $X\in\OD$. 
Then $X\ne\pu$, hence, there is 
a real $x\in X$. 
It follows from \osm\ that there is 
an ordinal $\la<\omi=\Om$, 
an element $f\in\koh_\la$, 
and an $\od$ map $H:\la^\om\to\bn$, 
such that $x=H(f)$. 
The set $P=\ens{f'\in\koh_\la}{H(f')\in X}$ is then \od\ 
and non-empty (contains $f$), and hence so is its image 
$Y=\ens{H(f')}{f'\in P}\sq X$ (contains $x$). 
Finally, $Y\in\odi$ by Lemma~\ref{L*}. 

To prove the second claim, let $X\ti Y$ be a condition in 
$\spe.$ 
By Lemma~\ref{sm6} there is a stronger saturated 
subcondition $X'\ti Y'\sq X\ti Y$.
By the first part of the lemma, let $X''\sq X'$ be a 
condition in $\odi$, and $Y''=Y'\cap\eke{X''}$.
Similarly, let $Y'''\sq Y''$ be a 
condition in $\odi$, and
$X'''=X''\cap\eke{Y'''}$.
Then $X'''\ti Y'''$ belongs to $\spei$.
\epf

\bcor
[assuming \osm]
\lam{egen}
If\/ $X\in\OD$ then there exists a \pge\ set\/ $G\sq\OD$ 
containing\/ $X$. 
If\/ $X\ti Y$ is a condition in $\spe$ then there exists 
a \ege\ set\/ $G\sq\spe$ containing\/ $X\ti Y$. 
\ecor
\bpf
By Lemma~\ref{den}, assume that $X\in\odi$. 
Then the set $\OD_{\sq X}$ of stronger conditions 
contains only countably many \od\ subsets by 
Lemma~\ref{2co}.
\epf

\subsection{The $\od$ forcing relation}
\las{odfr}
        
The forcing notion $\OD$ will play the same role below as 
the Gandy -- Harrington forcing in \cite{hms,kanHMS}.
There is a notable technical difference: under \osm, 
\dodg\ sets exist in the ground Solovay-model universe 
 by Corollary~\ref{egen}. 
Another notable difference is connected with the forcing 
relation. 

\bdf
[assuming \osm]
\lam{frd}
Let $\vpi(x)$ be an \dd\Ord\rit{formula}, that is, a 
formula with ordinals as parameters.

A  condition $X\in \OD$ is said to \rit{\pfo} $\vpi(\dox)$
iff $\vpi(x)$ is true 
(in the Solovay-model set universe considered)
for any \pge\ real $x$. 

If $\rE$ is an\/ $\od$ \eqr\ on $\bn$ then 
a condition $X\ti Y$ in $\spe$ is said to 
\rit{\efo} $\vpi(\doxl,\doxr)$
iff $\vpi(x,y)$ is true 
for any \ege\ pair $\ang{x,y}$.\qed
\edf

\ble
[assuming \osm]
\lam{frl}
Given an\/ \dd\Ord formula\/ $\vpi(x)$ and a\/ 
\pge\ real\/ $x$, if\/ $\vpi(x)$ is true\/ 
{\rm (in the Solovay-model set universe considered)} 
then there is a condition $X\in \OD$ containing\/ $x$, 
which\/ \pfo s\/ $\vpi(\dox)$.

Let\/ $\rE$ be an\/ $\od$ \eqr\ on $\bn.$ 
Given an\/ \dd\Ord formula\/ $\vpi(x,y)$ and a\/ 
\ege\ pair\/ $\ang{x,y}$, if\/ $\vpi(x,y)$ is true 
then there is a condition in\/ $\spe$ 
containing\/ $\ang{x,y}$, 
which\/ \efo s\/ $\vpi(\doxl,\doxr)$.
\ele
\bpf
To prove the first claim, put $X=\ens{x'\in\bn}{\vpi(x')}$. 
But this argument does not work for $\spe$. 
To fix the problem, we propose a longer argument which 
equally works in both cases --- but we present it in the 
case of $\OD$ which is slightly simpler.

Formally the forcing notion $\OD$ does not belong to $\rL$. 
But it is order-isomorphic to a certain forcing notion 
$P\in\rL$, namely, the set $P$ of \rit{codes}\snos
{A code of an $\od$ set $X$ is a finite sequence of logical 
symbols and ordinals which correspond to a definition in the 
form $X=\ens{x\in\rV_\al}{\rV_\al\models\vpi(x)}$.} 
of $\od$ sets in $\OD$.
The order between the codes in $P$, 
which reflects the relation $\sq$ 
between the $\od$ sets themselves, is expressible in $\rL$, 
too. 
Furthermore dense $\od$ sets in $\OD$ correspond 
to dense sets in the coded forcing $P$ in $\rL$.

Now, let $x$ be \pge\ and $\vpi(x)$ be true. 
It is a known property of the Solovay model that there is 
another \dd\Ord formula\/ $\psi(x)$ such that 
$\vpi(x)$ iff $\rL[x]\models\psi(x)$. 
Let $g\sq P$ be the set of all codes of conditions $X\in\OD$ 
such that $x\in X$. 
Then $g$ is a \dd Pgeneric set over $\rL$ by the choice of $x$, 
and $x$ is the corresponding generic object. 
Therefore there is a condition $p\in g$ which \dd Pforces 
$\psi(\dox)$ over $\rL$. 
Let $X\in\OD$ be the \od\ set coded by $p$, so that $x\in X$. 
To prove that $X$ \odf s $\vpi(\dox)$, let $x'\in X$ be a 
\pge\ real.
Let $g'\sq P$ be the \dd Pgeneric set of all codes of 
conditions $Y\in\OD$ such that $x'\in Y$. 
Then $p\in g'$, hence
$\psi(x')$ holds in $\rL[x']$, by the choice 
of $p$. 
Then $\vpi(x')$ holds 
(in the Solovay-model set universe) by the choice of $\psi$, 
as required.
\epf

\bcor
[assuming \osm]
\lam{frc}
Given an\/ \dd\Ord formula\/ $\vpi(x)$, if\/ $X\in\OD$ 
does not\/ \pfo\/ $\vpi(\dox)$ then there is a condition 
$Y\in \OD\yd Y\sq X$, 
which\/ \pfo s\/ $\neg\:\vpi(\dox)$.
The same for\/ $\spe$.\qed 
\ecor

\subsection{Some similar and derived forcing notions}
\las{simd}

Some forcing notions similar to $\OD$ and $\spe$ will be 
considered:
\ben
\cenu
\itla{sif1}\msur
$\odw=\ens{Q\sq W}{\pu\ne Q\in\od}$, where $W\sq\bn$ or 
$W\sq\bn\ti\bn$ is an $\od$ set.
Especially, in the case when $W\sq{\rE}$, where 
$\rE$ is an \od\ \eqr\ on $\bn$ 
(that is, ${\ang{x,y}\in W}\imp{x\rE y}$) --- note that 
$\eke{\dom W}=\eke{\ran W}$ in this case.

\itla{sif3}\msur
$(\spe)_{\sq X\ti Y} = 
\ens{X'\ti Y'\in\spe}{X'\sq X\land Y'\sq Y}$, where 
$\rE$ is an \od\ \eqr\ on $\bn$ and $X\ti Y\in \spe$.

\itla{sif4}\msur
$\odw\ti_{\rE} \odsw X = 
\ens{P\ti Y}
{P\in \odw\land Y\in \odsw X\land \eke Y\cap\eke{\dom P}\ne\pu}$, 
where 
$\rE$ is an \od\ \eqr\ on $\bn,$ $W\sq{\rE}$ is $\od$, 
$X\in \OD$, and $\eke X\cap\eke{\dom W}\ne\pu$ 
(equivalently, $\eke X\cap\eke{\ran W}\ne\pu$).

\itla{sif5}\msur
$\odw\ti_{\rE} \odw = 
\ens{P\ti Q}
{P,Q\in \odw\land \eke {\dom P}\cap\eke{\dom Q}\ne\pu}$, 
where 
$\rE$ is an \od\ \eqr\ on $\bn$ and $W\sq{\rE}$ is $\od$.
\een 
They have the same basic properties as 
$\OD$ --- the forcing notions of the form \ref{sif1},  
or as $\spe$ 
--- \ref{sif3}, \ref{sif4}, \ref{sif5}.
This includes such results and concepts as 
\ref{genx}, \ref{den}, \ref{egen}, the 
associated forcing relation as in \ref{frd}, and 
\ref{frl}, \ref{frc}, with suitable and rather transparent 
corrections, of course.

\subsection{Compression lemma}
\las{CL}

A set $A\sq2\Lom$ is an antichain if its elements are pairwise 
\dd\su incomparable, that is, no sequence in $A$ properly extends 
another sequence in $A$. 
Clearly any antichain is linearly ordered by $\lexe$.

Let $\vT=\Om^+$; 
the cardinal successor of $\Om$ in both $\bL$, the ground model, 
and its \dd{\coll\Om}generic extension postulated by \osm\ to 
be the set universe; in the latter, $\Om=\omi$ and $\vT=\om_2$.

\ble
[compression lemma]
\lam{apal31}
Assume that\/ $\Om\le\vt\le\vT$ and\/ $X\sq2^\vT$ is the image of\/ 
$\bn$ via an\/ $\od$ map. 
Then there is an\/ $\od$ antichain\/ $A(X)\sq2\Lom$ and an\/ 
$\od$ isomorphism\/ $f:\stk{X}{\lexe}\onto \stk{A(X)}{\lexe}$.
\ele
\bpf
If $\vt=\vT$ then, 
as $\card X\le\card{\bn}=\Om$, there is an ordinal $\vt<\vT$ 
such that $x\res\vt\ne y\res\vt$ whenever $x\ne y$ belong to $X$ 
--- this reduces the case $\vt=\vT$ to the case 
$\Om\le\vt<\vT$. 
We prove the latter by induction on $\vt$.

The nontrivial step is the step $\cof\la=\Om$, so that let 
$\vt=\bigcup_{\al<\Om}\vt_\al,$ for an increasing 
$\od$ sequence of ordinals $\vt_\al.$ 
Let $I_\al=[\vt_\al,\vt_{\al+1}).$ 
Then, by the induction hypothesis, for any $\al<\Om$ the set 
$X_\al=\ans{S\res I_\al:S\in X}\sq 2^{I_\al}$ is 
\dd\msl order-isomorphic to an antichain\/ 
$A_\al\sq 2\Lom$ via an\/ $\od$ isomorphism\/ $i_\al,$ 
and the map, which sends\/ $\al$ to\/ $A_\al$ and\/ $i_\al,$ 
is\/ $\od$. 
It follows that the map, which sends each $S\in X$ to the 
concatenation of all sequences $i_\al(x\res I_\al)$, 
is an $\od$ \dd\msl order-isomorphism $X$ onto an antichain 
in $2^\Om.$ 
Therefore, in fact it suffices to prove the lemma in the 
case $\vt=\Om.$ 
Thus let $X\sq 2^\Om.$ 

First of all, note that each sequence $S\in X$ is $\rod$.  
Lemma 7 in \cite{ksol}  shows that, in this case, we have 
$S\in\rL[S\res\eta]$ for an ordinal $\eta<\Om.$ 
Let $\eta(S)$ be the least such an ordinal, and 
$h(S)=S\res{\eta(S)},$ so that $h(S)$ is a countable initial 
segment of $S$ and $S\in\rL[h(S)].$ 
Note that $h$ is still $\od$. 

Consider the set $U=\ran h=\ens{h(S)}{S\in X}\sq 2\Lom.$ 
We can assume that every sequence $u\in U$ has a limit length. 
Then $U=\bigcup_{\ga<\Om}U_\ga,$ where 
$U_\ga=U\cap 2^{\om\ga}$ 
($\om\ga$ is the the \dd\ga th limit ordinal). 
For $u\in U_\ga,$ let $\ga_u=\ga$. 

If $u\in U$ then by construction the set $X_u=\ans{S\in X:h(S)=u}$ 
is $\od(u)$ and satisfies $X_u\sq\rL[u]$. 
Therefore, it follows from the known properties of 
the Solovay model that $X_u$ belongs to $\rL[u]$ and 
is of cardinality\/ $\le\Om$ in\/ $\rL[u]$.
Fix an enumeration $X_u=\ans{S_u(\al):\ga_u\le\al<\Om}$ 
for all $u\in U$. 
We can assume that the map $\al,u\longmapsto S_u(\al)$ 
is $\od$. 

If $u\in U$ and $\ga_u\le\al<\Om$, then we define  
a shorter sequence, $s_u(\al)\in 3^{\om\al+1}$, as follows.
\ben
\renu 
\itla{w1}\msur
$s_u(\al)(\xi+1)=S_u(\al)(\xi)$ for any $\xi<\om\al$. 

\itla{w2}\msur
$s_u(\al)(\om\al)=1$.

\itla{w4}
Let $\da<\al.$ 
If $S_u(\al)\res\om\da=S_{v}(\da)\res\om\da$ for some $v\in U$ 
(equal to or different from $u$) then 
$s_u(\al)(\om\da)=0$ whenever $S_u(\al)\lex S_{v}(\da),$ and 
$s_u(\al)(\om\da)=2$ whenever $S_{v}(\da)\msl S_u(\al).$ 

\itla{w3}
Otherwise (\ie, if there is no such $v$),  
$s_u(\al)(\om\da)=1$.
\een
To demonstrate that \ref{w4} is consistent, we show that 
$S_{u'}(\da)\res\om\da=S_{u''}(\da)\res\om\da$ implies $u'=u''.$ 
Indeed, as by definition $u'\subset S_{u'}(\da)$ and  
$u''\subset S_{u''}(\da),$ $u'$ and $u''$ must be 
\dd\sq compatible: let, say, $u'\sq u''.$ 
Now, by definition, $S_{u''}(\da)\in\rL[u''],$ therefore 
$\in \rL[S_{u'}(\da)]$ because 
$u''\sq S_{u''}(\da)\res\om\da=S_{u'}(\da)\res\om\da,$ 
finally $\in \rL[u'],$ which shows that  
$u'=u''$ as $S_{u''}(\da)\in X_{u''}$.

We are going to prove that the map 
$S_u(\al)\longmapsto s_u(\al)$ 
is a \dd\lex order isomorphism, so that 
$S_v(\ba)\lex S_u(\al)$ implies $s_v(\ba)\lex s_u(\al)$. 

We first observe that $s_v(\ba)$ and $s_u(\al)$ are 
\dd\sq incomparable. 
Indeed assume that $\ba<\al.$ 
If $S_u(\al)\res\om\ba\ne S_v(\ba)\res\om\ba$ then clearly 
$s_v(\ba)\not\sq s_u(\al)$ by \ref{w1}. 
If $S_u(\al)\res\om\ba= S_v(\ba)\res\om\ba$ then 
$s_u(\al)(\om\ba)=0\text{ or }2$ by \ref{w4} 
while $s_v(\ba)(\om\ba)=1$ by \ref{w2}. 
Thus all $s_u(\al)$ are mutually \dd\sq incomparable, so that 
it suffices to show that conversely 
$s_v(\ba)\msl s_u(\al)$ implies $S_v(\ba)\msl S_u(\al)$. 
Let $\za$ be the least ordinal such that 
$s_v(\ba)(\za)< s_u(\al)(\za);$ then 
$s_u(\al)\res\za= s_v(\ba)\res\za$ 
and $\za\le\tmin\ans{\om\al,\om\ba}.$ 

The case when $\za=\xi+1$ is clear: then by definition 
$S_u(\al)\res\xi= S_v(\ba)\res\xi$ while 
$S_v(\ba)(\xi)< S_u(\al)(\xi),$ 
so let us suppose that $\za=\om\da,$ where 
$\da\le\tmin\ans{\al,\ba}.$ 
Then obviously $S_u(\al)\res\om\da= S_v(\ba)\res\om\da.$ 
Assume that one of the ordinals $\al,\,\ba$ is equal to 
$\da,$ say, $\ba=\da.$ 
Then $s_v(\ba)(\om\da)=1$ while $s_u(\al)(\om\da)$ is 
computed by \ref{w4}. 
Now, as $s_v(\ba)(\om\da)< s_u(\al)(\om\da)$, we conclude that 
$s_u(\al)(\om\da)=2,$ hence $S_v(\ba)\msl S_u(\al),$ as required. 
Assume now that $\da<\tmin\ans{\al,\ba}.$ 
Then easily $\al$ and $\ba$ appear in one and the same class 
\ref{w4} or \ref{w3} with respect to the $\da$. 
However this cannot be \ref{w3} because 
$s_v(\ba)(\om\da)\ne s_u(\al)(\om\da).$ 
Hence we are in \ref{w4}, so that, for some (unique) 
$w\in U$. 
$0=S_v(\ba)\msl S_w(\da)\lex S_u(\al)=2,$ as required. 

This ends the proof of the lemma, except for the fact 
that the sequences $s_u(\al)$ belong to $3\Lom,$ but improvement 
to $2\Lom$ is easy.
\epf

\subsection{The dichotomy}
\las{Bsm}

Here we begin the proof of Theorem~\ref{mt}.
{\ubf We assume \osm\ in the course of the proof}. 
And we assume that the ordering $\cle$ of the theorem is 
just $\od$ --- then so is the associated \eqr\ $\apr$ 
and strict order $\cl$.

Let $\cF$ be the set of all $\od$ 
\lr\ order preserving maps $F:{\stk{\bn}\cle}\to\stk{A}\lexe$, 
where $A\sq2\Lom$ is an $\od$ antichain.
Let
\bce
$x\rE y\;$ iff $\;\kaz F\in\cF\:(F(x)=F(y))$
\ece
for $x,y\in\bn$.
Then $\rE$ is an $\od$ \eqr, \dd\od smooth in the sense that 
it admits an obvious $\od$ reduction to the equality on the 
set $2^\cF$. 

\ble
\lam{Lef2}
If\/ $R(x,y)$ is an\/ $\od$ relation and\/ 
$\kaz x,y\:({x\rE y}\imp R(x,y))$ 
then there is a function\/ $F\in\cF$ such that\/ 
$\kaz x,y\:({F(x)=F(y)}\imp R(x,y))$. 
\ele
\bpf 
Clearly $\card\cF=\Theta=\Om^+$ and $\cF$ admits 
an $\od$ enumeration $\cF=\ens{F_\xi}{\xi<\vT}$.
If $x\in\bn$ then let 
$ 
f(x)=F_0(x)\we F_1(x)\we \dots\we F_\xi(x)\we\;\dots
$  
--- the 
concatenation of all sequences $F_\xi(x)$. 
Then $f:\stk{\bn}\cle\to\stk{X}\lexe$ is an \od\ LR order 
preserving map, where $X=\ran f=\ens{f(r)}{r\in\bn}\sq2^\vT$, 
and ${f(x)=f(y)}\imp R(x,y)$ by the construction. 
By Lemma~\ref{apal31} there is an \od\ isomorphism 
$g:\stk{X}\lexe\onto\stk{A}\lexe$ onto an antichain 
$A\sq2\Lom.$   
The superposition $F(x)=g(f(x))$ proves the lemma.
\epf


\ble
\lam{i1sm}
Let \od\ sets\/ $\pu\ne X,Y\sq\bn$ satisfy\/ $\eke X=\eke Y$. 
Then the set\/ 
$B=\ens{\ang{x,y}\in X\ti Y}{x\rE y\cj x\cle y}$ is non-empty,
$\dom B=X$, $\ran B=Y$.
\ele
\bpf
It suffices to establish $B\ne\pu$.
The \od\ set 
$$
X' = \ens{x'\in\bn}
{\sus x \in X\,({x'\rE x}\land x'\cle x)}
$$ 
is downwards \dd\cle closed in each \dd\rE class, 
and if $B=\pu$ then $X'\cap Y=\pu$. 
By Lemma~\ref{Lef2}, 
there is a function $F\in\cF$ such that 
$x\in X'\imp x'\in X'$ holds whenever $F(x)=F(x')$ and 
$x'\cle x$.
It follows that  the derived  function
$$
G(x)=\left\{
\bay{rcl}
F(x)\we 0\,,&\text{whenewer}& x\in X'\\[1ex]
F(x)\we 1\,,&\text{whenewer}& x\in\bn\bez X' 
\eay
\right.
$$
belongs to $\cF.$ 
Thus if $x\in X\sq X'$ and $y\in Y\sq\bn\bez X'$ then 
$G(x)\ne G(y)$ and hence $x\nE y$. 
In other words, $\eke X\cap\eke Y=\pu$, a contradiction.
\epf

We'll make use of the 
\od-forcing notions $\OD$ and $\spe$.

\ble
\lam{SMsmu}
Condition\/ $\bn\ti\bn$ \efo s\/ $\doxl\rE\doxr$.  
\ele
\bpf
Otherwise, by Lemma \ref{frl}, 
there is a function $F\in\cF$ and a condition 
$X\ti Y$ in $\spe$ which \efo s  
$F(\doxl)(\xi)=0\ne 1= F(\doxr)(\xi)$ for a certain   
ordinal $\xi<\Om$. 
We may assume that $X\ti Y$ is a saturated condition. 
Then easily $F(x)(\xi)=0 \ne 1= F(y)(\xi)$ holds for any pair 
$\ang{x,y}\in X\ti Y$, so that we have $F(x)\ne F(y)$ and 
$x\nE y$ 
whenever $\ang{x,y}\in X\ti Y$, which contradicts the choice 
of $X\ti Y$ in $\spe$.
\epf

{\ubf Case 1:} 
$\approx$ and $\rE$ coincide on  $\bn$, so that 
${x\rE y}\eqv {x\apr y}$ for $x,y\in\bn$. 
By Lemma~\ref{Lef2} there is a single function $F\in\cF$ 
such that $F(x)=F(y)$ implies $x\apr y$ for all $x,y\in\Ua$, 
as required for \ref{mt1} of Theorem~\ref{mt}.     
\vtm

{\ubf Case 2:} 
$\apr$ is a \rit{proper} subrelation of $\rE$, hence, 
the $\od$ set 
$$
\Uo=\ens{x\in\bn}{\sus y\in \bn\:(x\napr y\land x\rE y)}
$$ 
(the domain of singularity) is non-empty. 
It follows that $\Uo\in\OD$ and $\Uo\ti\Uo$ is a condition 
in $\spe$. 
We'll work towards \ref{mt2} of Theorem~\ref{mt}.


\parf{The domain of singularity}
\las{sing}

Since the set $\Uo$ belongs to $\OD$, there is a set 
$\Ua\in\odi$, $\Ua\sq\Uo$.
Then obviously $\Ua\ti\Ua$ belongs to $\spei$.

\ble
\lam{31sm}
Condition\/ $\Ua\ti\Ua$ 
\efo s that the reals\/ $\doxl$ and\/ $\doxr$ are\/ 
\dd\cle incomparable. 
\ele
\bpf
Suppose to the contrary that, by Corollary~\ref{frc}, 
a subcondition $X\ti Y$ in $\spe$ 
either \efo s $\doxl\apr\doxr$ or \efo s $\doxl\cl\doxr$. 
We will get a contradiction in both cases. 
Note that $X,Y\sq\Ua$ are non-empty $\od$ sets and 
$\eke X\cap\eke Y\ne\pu$.

\bcl
\lam{ppp}
The set\/ 
$\ppw=\ens{\ang{x,x'}\in X\ti X}{x\rE x'\land x'\napr x}$ 
is non-empty.
\ecl
\bpf
Suppose to the contrary that $\ppw=\pu$, 
so ${\rE}$ coincides with ${\apr}$ on $X.$
As $X\sq \Ua$, at least one of the $\od$ sets 
$$
Z =\ens{z}{\sus x\in X\,(z\rE x\land z\ncle x)}\,,\; 
Z'=\ens{z}{\sus x\in X\,(z\rE x\land x\ncle z)}
$$ 
is non-empty; assume that, say, $Z\ne\pu$. 
Consider the $\od$ set 
$$
U=\ens{z}{\sus x\in X\,(z\rE x\land z\cle x)}\,.
$$
Then $X\sq U$ and $U\cap Z=\pu$, 
$U$ is downwards \dd\cle closed while $Z$ is upwards 
\dd\cle closed in each \dd\rE class, therefore $y\ncle x$ 
whenever $x\in U\land y\in Z\land x\rE y$, and hence we have
$\eke U\cap\eke Z=\pu$ be Lemma~\ref{i1sm}.
Yet by definition $\eke X\cap\eke Z\ne\pu$ and $X\sq U$, 
which is a contradiction.
\epF{Claim}

Suppose that 
condition $X\ti Y$ \efo s $\doxl\apr\doxr$. 
As $\ppw\ne\pu$ by Claim~\ref{ppp},
the forcing $\odw$ of all non-empty $\od$ 
sets $P\sq \ppw$ adds pairs $\ang{x,x'}\in \ppw$ of \pge\ 
(separately) 
reals $x,x'\in X$ which satisfy 
$x'\rE x$ and $x'\napr x$. 
If $P\in\odw$ then obviously $\eke{\dom P}=\eke{\ran P}$. 
Consider a more complex forcing 
$\cP=\odw\ti_{\rE}\OD$ of all pairs $P\ti Y'$, 
where $P\in\odw$,  
$Y'\in\OD$, $Y'\sq Y$, and $\eke{\dom P}\cap\eke{Y'}\ne\pu$. 
For instance, $\ppw\ti Y\in\odw\ti_{\rE}\OD$.  
Then $\cP$ adds a pair $\ang{\doxl,\doxr}\in \ppw$ and another 
real $\dox\in Y$ such that both pairs $\ang{\doxl,\dox}$ and 
$\ang{\doxr,\dox}$ belong to $X\ti Y$ and are \ege, 
hence, we have 
$\doxl\apr\dox\apr\doxr$ by the choice of $X\ti Y$. 
On the other hand, $\doxl\napr\doxr$ since the pair belongs 
to $\ppw$, which is a contradiction.

Now suppose  that condition $X\ti Y$ \efo s $\doxl\cl\doxr$. 
The set
$$
B=\ens{\ang{x,y}\in X\ti Y}{{y\rE x}\land {y\cle x}}
$$ 
is non-empty by Lemma~\ref{i1sm}.
Consider the forcing $\odwx$ of all non-empty $\od$ sets 
$P\sq B$; 
if $P\in\odwx$ then obviously $\eke{\dom P}=\eke{\ran P}$. 
Consider a more complex forcing $\odwex$ of all products 
$P\ti Q$, where $P,Q\in\odwx$ and 
$
\eke{\dom P}\cap\eke{\dom Q}\ne\pu\,.
$
In particular $B\ti B\in \odwex$.

Let $\ang{x,y;x',y'}$ be a \dd\odwex generic quadruple in 
$B\ti B$, 
so that both $\ang{x,y}\in B$ and $\ang{x',y'}\in B$ are 
\dd\odwx generic pairs in $B$, 
and both $y\cle x$ and  $y'\cle x'$ hold by the definition 
of $B$.
On the other hand, an easy argument shows that both criss-cross 
pairs $\ang{x,y'}\in X\ti Y$ and $\ang{x',y}\in X\ti Y$ are 
\dd{\spe}generic, hence $x\cl y'$ and $x'\cl y$ by the 
choice of $X\ti Y$. 
Altogether  $y\cle x\cl y'\cle x'\cl y$, which is a 
contradiction.
\epf

\parf{The splitting construction}
\las{spl}

Our aim is to define, in the universe of \osm, 
a splitting system of sets which 
leads to a function $F$ satisfying \ref{mt2} of Theorem~\ref{mt}. 
Let
\bce
$\ppl=\ens{\ang{x,y}\in\Ua\ti\Ua}{ x\rE y\cj x\cle y}$;\quad 
$\ppl\ne\pu$ by Lemma~\ref{i1sm}.
\ece
The construction will involve three forcing notions:
$\OD$, $\spe$, and 
\bce
$\odwb$,\,
the collection of all non-empty $\od$ sets 
$P\sq \ppl$.
\ece 
We also consider the dense (by Lemma~\ref{den}) subforcings 
$\odi\sq\OD$, $\spei\sq\spe$ (see Section~\ref{smod}),
and 
$$
\odwbi= \ens{Q\in\odwb}{\text{$Q$ is \odk}}\sq\odwb\,.
$$
Now note the following. 
\ben
\item 
As ${\Ua}\in {\odi}$, the set $\cD$ of all 
sets {open dense} in the restricted forcing
$\OD_{\sq\Ua}$,
is countable  by Lemma~\ref{den}; 
hence we can fix an enumeration 
$\cD=\ens{D_n}{n\in\om}$ such that 
$D_n\sq D_m$ whenever $m<n$. 

\item 
As ${\Ua\ti\Ua}\in {\spei}$, the set $\cD'$ of all 
sets, open dense in the restricted forcing
$(\spe)_{\sq \Ua\ti\Ua}$, is countable as above; 
fix an enumeration 
$\cD'=\ens{D'_n}{n\in\om}$ s.\,t.\ 
$D'_n\sq D'_m$ for $m<n$. 

\item 
If ${Q}\in {\odwbi}$ then the set $\cD(Q)$ of all 
sets {open dense} in the restricted forcing
$\OD_{\sq Q}$,
is countable  by Lemma~\ref{den}; 
hence we can fix an enumeration 
$\cD(Q)=\ens{D_n(Q)}{n\in\om}$ such that 
$D_n(Q)\sq D_m(Q)$ whenever $m<n$. 
\een
The chosen enumerations are not necessarily \od, of course.

A pair $\ang{u,v}$ of strings $u,\,v\in 2^n$ is called 
{\it crucial\/} iff $u=1^k\we 0\we w$ and $v=0^k\we 1\we w$ for 
some $k<n$ and $w\in 2^{n-k-1}.$ 
Note that each pair of the form $\ang{1^k\we 0,0^k\we 1}$ is a 
minimal crucial pair, and if $\ang{u,v}$ is a crucial pair 
then so is $\ang{u\we i,v\we i}$, but not $\ang{u\we i,v\we j}$
whenever $i\ne j$.
The graph 
of all crucial pairs in $2^n$ is actually a chain connecting all 
members of $2^n.$ 

We are going to define, {\ubf in the assumption of \osm}, 
a system of sets $X_u\in\odi$, where 
$u\in 2\lom,$ and sets $Q_{uv}\in \odwbi$, $\ang{u,v}$ being a 
crucial pair in some $2^n,$ satisfying the following conditions:

\ben
\itemsep=1mm
\def\theenumi{(\arabic{enumi})}
\def\labelenumi{\theenumi}
\itla{xq1}\msur
$X_u\in\odi$ and $Q_{uv}\in \odwbi$;

\itla{xq2}\msur
$X_{u\we i}\sq X_u$;

\itla{xq3}\msur
$Q_{u\we i\,,\,v\we i}\sq Q_{uv}$;

\itla{xq4}
if $\ang{u,v}$ is a crucial pair in $2^n$  
then $\pri Q_{uv}=X_u$ and $\prt Q_{uv}=X_v$;

\itla{xq5}\msur
$X_u\in D_n$ whenever $u\in 2^{n+1}$;

\itla{xq6}
if $u\yi v\in 2^{n+1}$ and $u(n)\not=v(n)$ then 
$X_u\ti X_v\in D'_n$ and 
$X_u\cap X_v=\pu$.

\itla{xq7}
if $\ang{u,v}=\ang{1^k\we 0\we w,0^k\we 1\we w}$ 
is a crucial pair in $2^{n+1}$ and $k<n$ 
(so that $w$ in not the empty string) then   
$Q_{uv}\in D_n(Q_{1^k\we 0,0^k\we 1})$;
\een

\bre
\lam{eke}
It follows from \ref{xq4} that 
$\eke{X_u}=\eke{X_v}$ for all $u\yi v\in 2^n,$ 
because $Q_{uv}\sq\ppl\sq{\rE}$ and $u\yi v$ are 
connected in $2^n$ by a chain of crucial pairs.\qed   
\ere

{\ubf Why this implies the existence of a 
function as in \ref{mt2} of Theorem~\ref{mt}?}


First of all, if $a\in 2^\om$ then the sequence 
of sets $X_{a\res n}$ is \dd\OD generic by \ref{xq5}, 
therefore the intersection $\bigcap_{n\in\om}X_{a\res n}$ is 
a singleton by Proposition~\ref{genx}. 
Let $F(a)\in\bn$ be its only element. 

It does not take much effort to prove that $F$ is continuous 
and $1-1$. 

Consider any $a\yi b\in\dn$ satisfying $a\nEo b$. 
Then $a(n)\not=b(n)$ for infinitely many $n,$ hence the pair 
$\ang{F(a),F(b)}$ is \dd\spe generic by \ref{xq7}, 
thus $F(a)$ and $F(b)$ are \dd\cle incomparable by 
Lemma~\ref{31sm}.

Consider $a\yi b\in\dn$ satisfying $a<_0 b$. 
We may assume that 
$a$ and $b$ are \dd{<_0}neighbours, \ie, $a=1^k\we 0\we w$ while 
$b=0^k\we 1\we w$ for some $k\in\om$ and $w\in 2^\om.$ 
The sequence of sets $Q_{a\res n\,,\,b\res n},\msur$ 
$n>k,$ is \dd\odwb generic by 
\ref{xq6}, hence it results in a pair of reals 
satisfying $x\mek y.$ 
However $x=F(a)$ and $y=F(b)$ by \ref{xq4}.

\punk{The construction of a splitting system}
\las{css}

Now the goal is to define, {\ubf in the assumption of \osm}, 
a system of sets $X_u$ and $Q_{uv}$ satisfying  
\ref{xq1} -- \ref{xq7} above.
Suppose that the construction has been completed up to a level 
$n,$ and expand it to the next level. 
From now on $s,\,t$ will 
denote strings in $2^n$ while $u,\,v$ will denote strings 
in $2^{n+1}.$ \vom 

{\ubf Step 0\/}. 
To start with, we set $X_{s\we i}=X_s$ for all $s\in 2^n$ and 
$i=0,1,$ and $Q_{s\we i\,,\,t\we i}=Q_{st}$ whenever $i=0,1$ and 
$\ang{s,t}$ is a crucial pair in $2^n.$ 
For the initial crucial pair $\ang{1^n\we 0,0^n\we 1}$ at 
this level, let 
$Q_{1^n\we 0\,,\,0^n\we 1}=X_{1^n}\ti X_{0^n}$. 
The newly defined sets satisfy \ref{xq1} -- \ref{xq4} except 
for the requirement $Q_{uv}\in \odwbi$ in \ref{xq1} for the 
pair $\ang{u,v}= \ang{1^n\we 0,0^n\we 1}$.

This ends the definition of ``initial values'' of $X_u$ 
and $Q_{uv}$ at the 
\dd{(n\pone)}th level. 
The plan is to gradually shrink the 
sets in order to fulfill \ref{xq5} -- \ref{xq7}.\vom 

{\ubf Step 1\/}. 
We take care of item \ref{xq5}. 
Consider an arbitrary $u_0=s_0\we i\in 2^{n+1}.$ 
As $D_n$ is dense there is a set  
$X'\in D_n\yd X'\sq X_{u_0}.$ 
The intention is to take $X'$ 
as the ``new'' $X_{u_0}$. 
But this change has to be propagated 
through the chain of crucial pairs, in order to preserve \ref{xq4}. 

Thus put $X'_{u_0}=X'$. 
Suppose that $u\in 2^{n+1},$ 
a set $X'_u\sq X_u$ has been defined, and 
$\ang{u,v}$ is a crucial pair, $v\in 2^{n+1}$ being not yet 
encountered. 
Define $Q'_{uv}=(X'_u\ti\bn)\cap Q_{uv}$ and  
$X'_v=\prt Q'_{uv}$. 
Clearly \ref{xq4} holds for the ``new'' sets 
$X'_u\yi X'_v\yi Q'_{uv}$. 
Similarly if $\ang{v,u}$ is a crucial pair, then 
define $Q'_{vu}=(\bn\ti X'_u)\cap Q_{vu}$ and  
$X'_v=\dom Q'_{uv}$. 
Note that still $Q'_{1^n\we 0\,,\,0^n\we 1}=X'_{1^n}\ti X'_{0^n}$.

The construction describes how the original change from $X_{u_0}$ 
to $X'_{u_0}$ spreads through the chain of crucial pairs in 
$2^{n+1},$ resulting in a system of new sets, $X'_u$ and 
$Q'_{uv},$ which satisfy \ref{xq5} for the particular 
$u_0\in 2^{n+1}.$ 
We iterate this construction consecutively for all 
$u_0\in 2^{n+1},$ getting finally a system of sets satisfying 
\ref{xq5} (fully) and \ref{xq4}, which we  denote by $X_u$ and 
$Q_{uv}$ from now on.\vom

{\ubf Step 2\/}. 
We take care of item \ref{xq6}. 
Consider a pair of $u_0$ and $v_0$ 
in $2^{n+1},$ such that $u_0(n)=0$ and $v_0(n)=1$. 
By the density 
of $D'_n$, there is a set $X'_{u_0}\ti X'_{v_0}\in D'_n$  
included in $X_{u_0}\ti X_{v_0}$. 
We may assume that 
$X'_{u_0}\cap X'_{v_0}=\emps.$ 
(Indeed it easily follows from Claim~\ref{ppp} that 
there exist reals $x_0\in X_{u_0}$ and $y_0\in X_{v_0}$ satisfying 
$x_0\rE y_0$ but $x_0\ne y_0,$ say $x_0(k)=0$ while $y_0(k)=1$. 
Define 
\dm
X=\ans{x\in X_0:x(k)=0\cj\sus y\in Y_0\:(y(k)=1\cj x\rE y)}\,,
\dm 
and $Y$ correspondingly; then $\eke X=\eke Y$ and $X\cap Y=\pu$.)

Spread the change from $X_{u_0}$ to $X'_{u_0}$ and from $X_{v_0}$ 
to $X'_{v_0}$ through the chain of crucial pairs in $2^{n+1},$ by 
the method of Step 1, until the wave of spreading from $u_0$ 
meets 
the wave of spreading from $u_0$ at the crucial pair 
$\ang{1^n\we 0,0^n\we 1}$. 
This 
leads to a system of sets $X'_u$ and $Q'_{uv}$ which satisfy 
\ref{xq7} for the particular pair $\ang{u_0,v_0}$ and still 
satisfy \ref{xq6} possibly except for the crucial pair 
$\ang{1^n\we 0,0^n\we 1}$ (for which basically the set 
$Q'_{1^n\we 0\,,\,0^n\we 1}$ is not yet defined for this step). 

By construction the previous steps leave 
$Q_{1^n\we 0\,,\,0^n\we 1}$ 
in the form 
$X_{1^n\we 0}\ti X_{0^n\we 1}$, 
where $X_{1^n\we 0}$ and 
$X_{0^n\we 1}$ are the ``versions'' at the end of Step 1). 
We now have the new sets, $X'_{1^n\we 0}$ and 
$X'_{0^n\we 1},$ included in resp.\ $X_{1^n\we 0}$ and 
$X_{0^n\we 1}$ and satisfying 
$\eke{X'_{0^n\we 0}}=\eke{X'_{0^n\we 1}}$. 
(Indeed $\eke{X'_{u_0}}=\eke{X'_{v_0}}$ held 
at the beginning of the change.) 
Now we put 
$Q'_{1^n\we 0\,,\,0^n\we 1}=(X'_{1^n\we 0}\ti X'_{0^n\we 1})
\cap\ppl$. 
Then $Q'_{1^n\we 0\,,\,0^n\we 1}\in\odwb$, and we have 
$\dom Q'_{1^n\we 0\,,\,0^n\we 1}=X'_{1^n\we 0}$, 
$\ran Q'_{1^n\we 0\,,\,0^n\we 1}=X'_{0^n\we 1}$ 
by Remark~\ref{eke} and Lemma~\ref{i1sm}.

This ends the consideration of the pair 
$\ang{u_0,v_0}$.

Applying this construction consecutively for all pairs of 
$u_0$ and $v_0$ with $u_0(n)=0$, $v_0(n)=1$
(including the pair 
$\ang{1^n\we 0,0^n\we 1}$) we finally get a system of sets 
satisfying \ref{xq1} -- \ref{xq6}, except 
for the requirement $Q_{uv}\in \odwbi$ in \ref{xq1} for the 
pair $\ang{u,v}= \ang{1^n\we 0,0^n\we 1}$, --- 
and these sets will be 
denoted still by $X_u$ and $Q_{uv}$ from now on.\vom   

{\ubf Step 3\/}.
Now we take care of \ref{xq7}. 
Consider a crucial pair in $2^{n+1}$,
$$
\ang{u_0,v_0}=\ang{1^k\we0\we w,0^k\we 1\we w}\in 2^{n+1}.
$$ 
If $k<n$ then $\ang{u_0,v_0}\ne\ang{1^k\we 0,0^k\we 1}$, 
the set $Q_{1^k\we 0,0^k\we 1}\in \odwbi$ 
is defined at a previous level, and 
$Q_{u_0,v_0}\sq Q_{1^k\we 0,0^k\we 1}$.
By the density, there 
exists a set $Q'_{u_0,v_0}\in D_n(Q_{1^k\we 0,0^k\we 1})$,  
$Q'_{u_0,v_0}\sq Q_{u_0,v_0}$. 
If $k=n$ then $\ang{u_0,v_0}=\ang{1^n\we 0,0^n\we 1}$, and by 
Lemma~\ref{den} there is a set $Q'_{u_0,v_0}\in \odwbi$, 
$Q'_{u_0,v_0}\sq Q_{u_0,v_0}$.
                                       
In both cases define $X'_{u_0}=\pri Q'_{u_0,v_0}$ and 
$X'_{v_0}=\prt Q'_{u_0,v_0}$ and spread this change through the 
chain of crucial pairs in $2^{n+1},$ exactly as above.
Note that $\eke{X'_{u_0}}=\eke{X'_{v_0}}$ as sets in $\odwb$ are 
included in $\rE$. 
This keeps $\eke{X'_u}=\eke{X'_v}$ 
for all $u,\,v\in 2^{n+1}$ through the spreading.  
 
Executing this step for all crucial pairs in $2^{n+1},$ 
we finally accomplish the construction of a 
system of sets satisfying \ref{xq1} through \ref{xq7}.\vtm

\qeDD{Theorem~\ref{mt}}

{\small

}

\end{document}